\documentclass[12pt, twoside, leqno]{article}

\usepackage{amsmath,amsthm}
\usepackage{amssymb}
\usepackage{arydshln}

\usepackage{enumerate}

\usepackage{graphicx}

\newtheorem{thm}{Theorem}[section]
\newtheorem{cor}[thm]{Corollary}
\newtheorem{lem}[thm]{Lemma}
\newtheorem{prop}[thm]{Proposition}

\theoremstyle{definition}
\newtheorem{defn}[thm]{Definition}
\newtheorem{rem}[thm]{Remark}

\numberwithin{equation}{section}

\usepackage{pgf,tikz}
\usetikzlibrary{arrows}

\begin{document}

\baselineskip=12.1pt


\title{Some remarks on the square graph of the hypercube }

\author{S.Morteza Mirafzal\\
Department of Mathematics \\
  Lorestan University, Khorramabad, Iran\\
\\
E-mail: mirafzal.m@lu.ac.ir\\
E-mail: smortezamirafzal@yahoo.com}

\date{}

\maketitle


\renewcommand{\thefootnote}{}



\footnote{\emph{Date}:  }

\renewcommand{\thefootnote}{\arabic{footnote}}
\setcounter{footnote}{0}
\date{}
\begin{abstract}
Let $\Gamma=(V,E)$ be a graph. The square graph $\Gamma^2$ of the graph $\Gamma$  is the graph with the vertex set $V(\Gamma^2)=V$ in which
two vertices are adjacent if and only if their distance in $\Gamma$ is at most two. The square graph of the hypercube $Q_n$ has some interesting properties. For instance,  it is highly symmetric and panconnected.
 In this paper, we investigate some   algebraic properties of the graph ${Q^2_n}$.   In particular,   we show that the graph ${Q^2_n}$ is distance-transitive.     We show that the graph ${Q^2_n}$ is an imprimitive distance-transitive graph if and only if $n$ is an odd integer. Also, we determine the spectrum of the graph $Q_n^2$.  Finally, we show that when $n >2$ is an even integer, then ${Q^2_n}$ is an automorphic graph, that is, $Q_n^2$ is a distance-transitive primitive graph which is not a complete or a  line graph.
\end{abstract}
\maketitle
\section{ Introduction}
\noindent
considered as an undirected simple graph where $V=V(\Gamma)$ is the vertex-set
and $E=E(\Gamma)$ is the edge-set. For all the terminology and notation
not defined here, we follow \cite{1,3,5,6,9}.\

Let $\Gamma=(V,E)$ be a graph.  The $square$ $graph$   $\Gamma^{2}$ of the graph $\Gamma$  is the (simple)  graph with vertex set $V$ in which two vertices are adjacent if and only if their distance in $\Gamma$ is at most two. It is easy to see that $Aut(\Gamma) \leq Aut(\Gamma^2)$, where $Aut(\Gamma)$ denotes the automorphism group of the graph $\Gamma$. Thus, if the graph $\Gamma$ is a vertex-transitive graph, then $\Gamma^2$ is a vertex-transitive graph. A graph $\Gamma$ of order $n > 2$ is $Hamilton$-$connected$ if for
any pair of distinct vertices $u$ and $v$, there is a Hamilton $u$-$v$ path, namely, there is
a $u$-$v$ path of length $n-1$. It is clear that if a graph $\Gamma$ is Hamilton-connected  then it is Hamiltonian. A graph $\Gamma$ of order $n > 2$ is panconnected if for every two vertices $u$ and $v$, there is
a $u$-$v$ path of length $l$ for every integer $l$ with
$d(u,v) \leq l \leq n-1$. Note that if a graph $\Gamma$ is   panconnected, then it is Hamilton-connected.
 It is a well known fact that when a graph $\Gamma$
is 2-connected, then its square $\Gamma^2$ is panconnected   \cite{4,7}. Using this fact, and  an algebraic property of Johnson graphs, recently it has been proved
that the Johnson graphs are panconnected  \cite{17}.
\

Let $ n \geq 2 $ be an integer.
The hypercube  $Q_n$ is the graph whose vertex-set is $ \{0,1  \}^n $, where two $n$-tuples  are adjacent if  they differ in precisely one coordinate. This graph has  been studied from various   aspects by many authors. Some recent works concerning some algebraic aspects of this graph include \cite{13,15,19,28}. It is a well known fact that the graph $Q_n$ is a distance-transitive graph \cite{1,3}, and hence it is edge-transitive. Now, using a well known result due to Watkins \cite{27}, it follows that
the connectivity of $Q_n$ is maximal, that is, $n$.  Like the hypercube $Q_n$, its square, namely,  the graph $Q_n^2$ has some interesting properties. For instance, when $n\geq 2$, then $Q_n$ is 2-connected. Now using a known result due to  Chartrand and Fleischner \cite{4,7}, it follows that   $Q_n^2$ is a panconnected graph.
 Also, since $Q_n$ is vertex-transitive,   the graph $Q_n^2$ is vertex-transitive, as well. Hence $Q_n^2$ is a
regular graph  and it is easy to check that its valency is $n$+$n\choose 2$=$n+1\choose 2$. If
$n=2$, then $Q_n^2$ is the complete graph $K_4$.
When $n=3$, then $Q_n^2$ is a 6-regular graph with 8 vertices. This graph is isomorphic with a graph known as the $coktail$-$party$ graph $CP(4)$ \cite{1}. It can be shown that when $n=4$, then the graph   $Q_n^2$ is a 10-regular graph with 16 vertices, which  is isomorphic to the complement of  the graph known as the $Clebsch$ graph \cite{9}.

    In this paper, we determine the automorphism group of the graph $Q_n^2$. Then we show that $Q_n^2$ is a distance-transitive graph.
  This implies that the connectivity of the graph $Q_n^2$ is maximal, namely, its valency     $n+1\choose 2$.
Also,  we will see that the graph ${Q^2_n}$ is an imprimitive distance-transitive graph if and only if $n$ is an odd integer. A graph $\Gamma$ is called an $automorphic$ graph, when it is a distance-transitive primitive graph which is not a complete or a  line graph \cite{1}. In the last step of the paper, we show that the graph $Q_n^2$ is an automorphic graph if and only if $n$ is an even integer.

\section{Preliminaries}
The graphs $\Gamma_1 = (V_1,E_1)$ and $\Gamma_2 =
(V_2,E_2)$ are called $isomorphic$, if there is a bijection $\alpha
: V_1 \longrightarrow V_2 $   such that  $\{a,b\} \in E_1$ if and
only if $\{\alpha(a),\alpha(b)\} \in E_2$ for all $a,b \in V_1$.
In such a case the bijection $\alpha$ is called an $isomorphism$.
An $automorphism$ of a graph $\Gamma$ is an isomorphism of $\Gamma
$ with itself. The set of automorphisms of $\Gamma$  with the
operation of composition of functions is a group  called the
$automorphism\  group$ of $\Gamma$ and denoted by $ Aut(\Gamma)$.

 The
group of all permutations of a set $V$ is denoted by $Sym(V)$  or
just $Sym(n)$ when $|V| =n $. A $permutation$ $group$ $G$ on
$V$ is a subgroup of $Sym(V).$  In this case we say that $G$ $acts$
on $V$. If $G$ acts on $V$ we say that $G$ is
$transitive$ on $V$ (or $G$ acts $transitively$ on $V$) if given any two elements $u$ and $v$ of
$V$, there is an element $ \beta $ of  $G$ such that  $\beta (u)= v
$.  If $\Gamma$ is a graph with vertex-set $V$  then we can view
each automorphism of $\Gamma$ as a permutation on $V$  and so $Aut(\Gamma) = G$ is a
permutation group on $V$.

A graph $\Gamma$ is called $vertex$-$transitive$ if  $Aut(\Gamma)$
acts transitively on $V(\Gamma)$. We say that $\Gamma$ is $edge$-$transitive$ if the group $Aut(\Gamma)$ acts transitively  on the edge set $E$, namely, for any $\{x, y\} ,   \{v, w\} \in E(\Gamma)$, there is some $\pi$ in $Aut(\Gamma)$,  such that $\pi(\{x, y\}) = \{v, w\}$.  We say that $\Gamma$ is $symmetric$ (or $arc$-$transitive$) if  for all vertices $u, v, x, y$ of $\Gamma$ such that $u$ and $v$ are adjacent, and also, $x$ and $y$ are adjacent, there is an automorphism $\pi$ in $Aut(\Gamma)$ such that $\pi(u)=x$ and $\pi(v)=y$. We say that $\Gamma$ is $distance$-$transitive$ if  for all vertices $u, v, x, y$ of $\Gamma$ such that $d(u, v)=d(x, y)$, where $d(u, v)$ denotes the distance between the vertices $u$ and $v$  in $\Gamma$,  there is an automorphism $\pi$ in $Aut(\Gamma)$ such that  $\pi(u)=x$ and $\pi(v)=y.$\

A vertex cut of the graph $\Gamma$ is a subset $U$ of $V$ such that the subgraph
$\Gamma-U$ induced by the set $V-U$ is either trivial or not connected.  The $connectivity$
$\kappa(\Gamma)$ of a nontrivial connected graph $\Gamma$ is the minimum cardinality of all vertex
cuts   of $\Gamma$. If we denote by $\delta(\Gamma)$ the minimum degree of $\Gamma$, then $\kappa(\Gamma) \leq \delta(\Gamma)$. A graph $\Gamma$ is called $k$-$connected$ (for $k \in \mathbb{N}$) if $|V(\Gamma)| > k$ and $\Gamma-X$ is connected for every subset $X \subset V(\Gamma)  $ with $|X|< k$. It is trivial that if a  positive integer $m$ is such that $m \leq \kappa(\Gamma)$, then $\Gamma$ is an $m$-connected graph.   We have the following fact.

\begin{thm} $\cite{27}$ \label{a} If a connected graph $\Gamma$ is edge-transitive, then $\kappa(\Gamma) = \delta(\Gamma)$, where $\delta(\Gamma)$ is the minimum degree of vertices of $\Gamma$.

\end{thm}

Let  $n,k \in \mathbb{ N}$ with $ k < n,   $ and let $[n]=\{1,...,n\}$. The $Johnson\  graph$ $J(n,k)$ is defined as the graph whose vertex set is $V=\{v\mid v\subseteq [n], |v|=k\}$ and two vertices $v$,$w  $ are adjacent if and only if $|v\cap w|=k-1$.   The class of  Johnson  graphs is a well known class of   distance-transitive graphs  \cite{3}.   It is an easy task to show that the set   of mappings  $H= \{ f_\theta \ | \  \theta \in$ S$ym([n]) \} $,  $f_\theta (\{x_1, ..., x_k \}) = \{ \theta (x_1), ..., \theta (x_k) \} $,    is a subgroup of $ Aut( J(n,k) ) $ \cite{9}.   It has been shown that   $Aut(J(n,k)) \cong$ S$ym([n])$  if  $ n\neq 2k, $  and $Aut(J(n,k)) \cong$ S$ym([n]) \times \mathbb{Z}_2$, if $ n=2k$,   where $\mathbb{Z}_2$ is the cyclic group of order 2 \cite{3,12,23}.

Although
in most situations  it is difficult  to determine the automorphism group
of a graph $\Gamma$  and how it acts on its vertex and edge sets,    there are various papers in the literature, and some of the recent works
include \cite{8,10,12,13,14,16,18,19,20,21,22,24,26,28}.\

Let $G$ be any abstract finite group with identity $1$, and
suppose $\Omega$ is a subset of   $G$, with the
properties:

(i) $x\in \Omega \Longrightarrow x^{-1} \in \Omega$,   $ \ (ii)
 \ 1\notin \Omega $.

The $Cayley\  graph$  $\Gamma=Cay (G; \Omega )$ is the (simple)
graph whose vertex-set and edge-set  are defined as follows:

$V(\Gamma) = G $,  $  E(\Gamma)=\{\{g,h\}\mid g^{-1}h\in \Omega \}$.\\It can be shown that the
 Cayley   graph   $\Gamma=Cay (G; \Omega )$ is connected if and only if the set $\Omega$ is a generating set in the group $G$ \cite{1}.\

The group $G$ is called a semidirect product of $ N $ by $Q$,
denoted by $ G=N \rtimes Q $,
 if $G$ contains subgroups $ N $ and $ Q $ such that:  (i)$
N \unlhd G $ ($N$ is a normal subgroup of $G$); (ii) $ NQ = G $; and
(iii) $N \cap Q =1 $. \

\section{Main results}

The hypercube  $Q_n$ is the graph whose vertex set is $ \{0,1  \}^n $, where two $n$-tuples  are adjacent if  they differ in precisely one coordinates. It is easy  to show that  $ Q_n = Cay(\mathbb{Z}_{2}^n; S )$, where $\mathbb{Z}_{2}$ is
the cyclic group of order 2, and $S=\{ e_i \  | \  1\leq i \leq n \}, $ where  $e_i = (0, ..., 0, 1, 0, ..., 0)$,  with 1 at the $i$th position. It is   easy   to show that the set  $H= \{ f_\theta |   \theta \in Sym([n]) \} $, $ f_\theta (x_1, ..., x_n ) =( x_{\theta (1)}, ...,  x_{\theta (n)}) $  is a subgroup of the group $Aut(Q_n)$.  It is clear that $H \cong Sym([n])$. We know that in every Cayley graph $\Gamma= Cay(G;S)$, the group $Aut(\Gamma)$ contains a subgroup isomorphic with the group $G$.  In fact,  if $x\in \mathbb{Z}_2^n$, and we define the mapping $f_x(v)=x+v$, for every $v\in  V(Q_n)$, then $f_x$ is an automorphism of the hypercube $Q_n$.  Hence $\mathbb{Z}_{2}^n $ is (isomorphic with) a subgroup of $Aut(Q_n)$. It has been proved that $Aut(Q_n) = \langle \mathbb{Z}_{2}^n, Sym([n]) \rangle \cong \mathbb{Z}_{2}^n \rtimes Sym([n])$ \cite{13}. It is clear that when $\Gamma$ is a graph then $Aut(\Gamma)$ is a subgroup of $Aut(\Gamma^2)$.
Thus we have $Aut(Q_n)\leq Aut(Q_n^2)$. In the sequel, we wish to show that the graph $Q_n^2$ is a distance-transitive graph, and for doing this we need the automorphism group of $Q_n^2$. When $n=3$, then $Q_n^2$ is isomorphic with the coktail-party graph $CP(4)$. The complement of this graph is a disjoint union of 4 copies of $K_2$.  Thus $Aut(Q_3^2) \cong Sym([2])\  wr_I \   Sym([4]) $, where $I=\{1,2,3,4  \}$ \cite{3,21} (for an acquaintance with the notion of wreath product of groups see \cite{6}). Now it can be checked that this graph is a distance-transitive graph. Hence, in the sequel we assume that $n \geq 4$.
It is easy to see that for the   graph $Q_n^2$ we have,  $Q_n^2=Cay(\mathbb{Z}_2^n;T)$,   $T=S \cup S_1$, where $S_1= \{e_i+e_j | \  i,j \in [n], i\neq j  \}$.  Let $A=Aut(Q_n^2)$ and $A_0$ be the stabilizer subgroup of the vertex $v=0$ in $A$.  Since $ Q_n^2$ is a vertex-transitive graph, then from the orbit-stabilizer theorem we have $|A|=|A_0||V(Q_n^2)|=2^n|A_0|$. The following lemma determines an upper bound for $|A_0|$.
%
%
\begin{lem} \label{c1}
Let $n > 4$ and $A=Aut(Q_n^2)$. Let  $A_0$ be the stabilizer subgroup of the vertex $v=0$. Then $|A_0|\leq(n+1)!$.
\end{lem}

\begin{proof}   Let $\Gamma=Q_n^2$.  We know that  $  \Gamma=Cay(\mathbb{Z}_2^n;T)$,   $T=S \cup S_1$, where $S=\{e_i \ | \ 1 \leq i \leq n  \}$ and $S_1= \{e_i+e_j | \  i,j \in [n], i\neq j  \}$.  Let $f\in A_0$. Then $f(T)=T$. Let $G$ be the subgraph of $\Gamma$ which is induced by the subset $T$. Let $h=f|_T$ be the restriction of the mapping $f$ to the subset $T$. It is clear that $h$ is an  automorphism of the graph $G$.
It is easy to see that the mapping $\Phi:A_0 \rightarrow Aut(G)$, which is  defined by the rule $\Phi(g)=g|_T$,  is a group homomorphism.  Thus we have $\frac{A_0}{ker(\Phi)} \cong im(\Phi)$, and hence we have $|A_0|=|ker(\Phi)||im(\Phi)|$. Since $im(\Phi)$ is a subgroup of $Aut(G)$, then we have $|A_0|\leq |ker(\Phi)||Aut(G)|$. If we show that $|Aut(G)|\leq(n+1)!$ and  $ker(\Phi)=\{1\}$,  then  the lemma is proved. Hence in the rest of the proof  we show that:
\\(i) $|Aut(G)|\leq (n+1)!$, \\ (ii) $ker(\Phi)=\{1\}.$ \\
(i)  We give two proofs for proving this claim. The first is more elementary than the second, but  we need some parts of it  in the proof of (ii). The second is based on  the automorphism group of the Johnson graph $J(n,k)$.\\
{\bf Proof 1 of (i)}.  Consider the graph $G$. In $T=V(G)$, consider the subgraphs induced by the subsets $C_0=S=\{e_i |\ 1 \leq i \leq n  \}$, $C_i=\{e_i, e_i+e_j| \  1\leq j\leq n, i\neq j  \}$, $1 \leq i \leq n $ (we also denote by $C_i$ the subgraph induced by the set $C_i$ ).  It is clear that $C_0$   is an $n$-clique in the graph $G$.  Note that if $e_i+e_r$ and $e_i+e_s$ are two elements of $C_i$, then we have $(e_i+e_r)-(e_i+e_s)= e_r+e_s \in T$. Hence each $C_i$ is also an  $n$-clique in the graph $G$. It can be shown that each $C_i$, $0\leq i \leq n$ is a maximal $n$-clique in $G$. It is clear that if $i\neq 0$, then $C_0 \cap C_i =\{ e_i \}$. Moreover, if $i,j \in \{1,...,n\}$ and  $i\neq j$,  then  $C_i \cap C_j = \{e_i+e_j \}$. Let $M$ be a maximal $n$-clique in the graph $G$.
It is not hard to show that $M=C_j$ for some $j \in \{ 0,1, ..., n \}$. If $a$ is an automorphism of the graph $G$, then $a(C_j)$ is a maximal $n$-clique in the graph $G$. Hence the natural action of $a$ on the set $X=\{ C_0,C_1,...,C_n \}$ is a permutation on $X$. Let $G_1$ be the graph with the vertex set $X$ in which two vertices $v$ and $w$ are adjacent if and only if $v \cap w \neq \emptyset$. Now, it is clear that $G_1 \cong K_{n+1}$, the complete graph on $n+1$ vertices, and hence $Aut(G_1) \cong Sym(X)$. Let $a\in Aut(G)$ be such that
$a(C_j)=C_j$, for each $j \in \{0,1,...,n\}$. Noting that $C_0 \cap C_i=\{ e_i \}$, $i\neq 0$, we
deduce that $a(x)=x$ for every $x\in C_0$. Note that the vertex  $e_i+e_j$ is the unique common neighbor of vertices $e_i$ and $e_j$ in the graph $G$  which is not in $C_0$.  This implies  that
$a(e_i+ e_j) =e_i+e_j$. Therefore we have  $a(v)=v$ for every $v\in T$.   Now it is easy to see that the mapping $\pi: Aut(G)\rightarrow Aut(G_1) $ defined by the rule $\pi (a)=f_a$, where  $f_a(C_i)=a(C_i)$ for every $C_i\in X$,  is an injection and therefore we have $ (n+1)!\geq |Aut(G)| $. \

{\bf Proof 2 of (i)}. Consider the graph $G$. We
show that this graph is isomorphic with the  Johnson
graph $J(n+1,2)$. We define the mapping
$$f: V(G)\rightarrow V(J(n+1,2)), $$   by the rule:
$$ f(v) =  \begin{cases}
\{i,n+1\},     if \ v =e_i \\   \{ i,j \}, \ \ \  \ \      if\ \  v=e_i+e_j\\
 \end{cases} $$
It is clear that $f$ is a bijection. Let $\{v,w\}$ be an edge in the graph $G$. Then we have three possibilities: \\(1)
$\{v,w\}=\{e_i,e_j  \}$, (2) $\{v,w\}=\{e_i,e_i +e_k  \}$, (3) $\{v,w\}=\{ e_i+e_k,e_i +e_j \}$. \\ Now,
we have (1) $f(\{v,w\})=\{ \{i,n+1  \},\{j,n+1  \} \}$, (2) $f(\{v,w\})=\{ \{i,n+1  \},\{i,k\} \}$,  (3) $f(\{v,w\})=\{ \{i,k\},\{j,k\} \}$. It follows   that $f$ is a graph isomorphism. Hence, $Aut(G) \cong Aut(J(n+1,2))$. Since $Aut(J(n+1,2)) \cong Sym([n+1])$ \cite{3,12,23}, then we have $Aut(G) \cong Sym([n+1])$. \

(ii) we now show that $ker(\Phi)=\{ 1 \}$. Let $f \in ker(\Phi)$. Then $f(0)=0$ and  $h=f|_T$ is the identity automorphism of the graph $G$. Hence $f(x)=x$ for every $x \in T$. Note that when $x\in T$, then   $w(x) \in \{1,2 \}$, where $w(x)$ is the weight of $x$, that is, the number of 1s in the $n$-tuple $x$. Let $x\in V(\Gamma)$ and $w(x)= m$. We show by induction on $m$, that $f(x)=x$. It is clear that when $m=0,1,2$, then the claim is true. Let the claim be true when $w(x)\leq m$, $m\geq2$. We show that if $w(x)=m+1$, then $f(x)=x$. Let $y=e_{i_1}+...+e_{i_m}+e_{i_{m+1}}$ be a vertex of weight $m+1$. Let $v=y+e_{i_m}+e_{i_{m+1}}$. Since $W(v)=m-1$, thus $f(v)=v$. Let $N$ be the subgraph of $\Gamma$ which is induced by the set $N(v)$.  Since $\Gamma$ is vertex-transitive, then $G \cong N$.
Also,  since $f(v)=v$, then the restriction of $f$ to $N(v)$ is an automorphism of the graph $N$.
In $N(v)$ we define the subsets $M_0=\{ v+e_i |\
1\leq i \leq n \}$, $M_i=\{v+e_i, v+e_i+e_j |\  1\leq j \leq n , j \neq i\}$,  $1\leq i \leq n$. It can be check that the subgraph induced by each $M_i$
  is a maximal $n$-clique in the graph $N$. Also,  $M_0 \cap M_i=\{v+e_i \}$. Moreover,  $v+e_i+e_j$ is the unique common neighbor of the vertices $v+e_i$ and $v+e_j$ in the graph $N$ which is not in $M_0$.  If $x\in M_0$, then $f(x)=x$,
because $w(x)\leq m$. This implies that $f(M_i)=M_i$. Now,  by an argument similar to what is done in    proof 1, we can see  that $f(x)=x$ for every $x\in N(v)$. Since $y\in N(v)$, we have $f(y)=y$. We now conclude that $f$ is the identity automorphism of $\Gamma$. Hence we have $ker(\Phi)=\{ 1 \}$.
\end{proof}
%
%
\begin{thm} \label{c2}
Let $n > 4$ and $\Gamma =Q_n^2$ be the square of the hypercube $Q_n$. Then we have
$Aut(\Gamma) \cong  \mathbb{Z}_{2}^n \rtimes Sym([n+1])$.
\end{thm}

\begin{proof}
Let  $A_0$ be the stabilizer subgroup of the vertex $v=0$ in the group $Aut(\Gamma)$. We know from Lemma \ref{c1}, that $|A_0| \leq (n+1)!$.
Let  $T$ and $X=\{C_0,...,C_n\}$ be the sets which are defined in the proof of Lemma \ref{c1}. Note that $\mathbb{Z}_{2}^n$ is a  vector space  over the field $\mathbb{Z}_2$ and $C_i$, $0 \leq i \leq n$,  is a basis for this vector space.  Let $f_i : C_0 \rightarrow C_i$ be a bijection. We can linearly extend $f_i$ to an automorphism $e(f_i)$ of the group $\mathbb{Z}_{2}^n$. It is clear that $e(f_i) \in A_0$.  We know that every automorphism of the group $\mathbb{Z}_{2}^n$ which fixes the set $T$ is an automorphism of the graph $\Gamma$. We can see that when $x,y\in C_i$ and $x\neq y$ then $x+y \in T$. Thus we have $e(f_i)(e_r+e_s)=e(f_i)(e_r)+e(f_i)(e_s) \in T$. Hence we have $e(f_i)(T)=T$. Since the number of permutations $f_i$ is $n!$, hence the number of automorphisms of $e(f_i)$ is $n!$. Note that when $i \neq j$, then $e(f_i) \neq e(f_j)$. Now since $0\leq i \leq n$, then we have at least $(n+1)(n!)=(n+1)!$ distinct automorphisms in the group $A_0$. Thus by Lemma \ref{c1}, we have $|A_0|=(n+1)!$.
 We saw, in the proof of Lemma \ref{c1}, that $A_0$ is isomorphic with a subgroup of $Sym([n+1])$. Hence we deduce that $A_0 \cong Sym([n+1])$.\

 We know, by the orbit-stabilizer theorem, that   $|V(\Gamma)||A_0|=|Aut(\Gamma)|$. Thus
we have $|Aut(\Gamma)|=2^n[(n+1)!]$. For every
$v \in \mathbb{Z}_2^n$,  the mapping $ f_v(x)=v+x$,
  for every $x\in \mathbb{Z}_2^n$, is an
automorphism of the graph $\Gamma$. It is easy to check that $ L=\{f_v|\  v\in \mathbb{Z}_2^n  \}$
is a subgroup of $Aut(\Gamma)$ which is isomorphic with $\mathbb{Z}_2^n$. Also it is easy to check that $L \cap A_0 =\{ 1 \}$. Hence we have $|LA_0|=|L||A_0|=2^n[(n+1)!]=|Aut(\Gamma)|$. This implies that $Aut(\Gamma)=LA_0$.  Also we can
see that for every $v\in \mathbb{Z}_2^n$ and every $a\in A_0$  we have
$ a^{-1}f_va=f_{a^{-1}(v)}$. Thus we deduce that $L$ is a normal subgroup of $Aut(\Gamma)$. We now conclude that 
$$Aut(\Gamma) \cong L \rtimes A_0 \cong \mathbb{Z}_{2}^n \rtimes Sym([n+1]).$$
\end{proof}
The graph $Q_n^2$ has some interesting properties. In the next theorem, we show that $Q_n^2$ is distance-transitive.
%
%
\begin{thm} \label{c3} Let $n \geq 4$ be an integer. Then the graph $Q_n^2$ is a distance-transitive graph.

\end{thm}

\begin{proof}  Let $v$ and $w$ be vertices in $Q_n^2$. It is easy to check that $d(x,y)=\lceil \frac{w(x+y)}{2} \rceil$. Hence we have $d(x,0)=\lceil \frac{w(x)}{2}\rceil $. Let $D$ be the diameter of $Q_n^2$. it follows from the first two sentences that $D=\lceil\frac{n}{2} \rceil$. Let $A_0$ be the stabilizer subgroup   the vertex $v=0$ in $Aut(Q_n^2)$.    Since the graph $Q_n^2$ is a vertex-transitive graph,  it is sufficient to show that the action of $A_0$ on the set $\Gamma_k$  is transitive, where $\Gamma_k$ is the set of vertices at distance $k$ from the vertex $v=0$. Let $x$ and $y$ be two vertices in $\Gamma_k$. There are two possible cases, that is, (i) $w(x)=w(y)$ or  (ii) $w(x) \neq w(y)$.\\ (i) Let  $w(x)=w(y)$. We know that $w(x) \in \{2k,2k-1 \}$. Without loss of generality, we can assume that $w(x)=2k$.  Let $x=e_{i_1}+...+e_{i_{2k}}$  and  $y=e_{j_1}+...+e_{j_{2k}}$. There are vertices $e_{x_1},...,e_{x_{n-2k}}$ and
$e_{y_1},...,e_{y_{n-2k}}$ in $Q_n^2$ such that\\
$\{ e_{i_1},...,e_{i_{2k}},e_{x_1},...,e_{x_{n-2k}}  \}$=$C_0$=$ \{e_1,e_2,...,e_n \}$=$\{e_{j_1},...,e_{j_{2k}},e_{y_1},...,e_{y_{n-2k}}  \}$. Let $f$ be the permutation on the set $C_0$ which is defined by the rule, $f(e_{i_r})=e_{j_r}, 1 \leq r \leq 2k$, and $f(e_{x_{l}})=e_{y_{l}}, 1 \leq l \leq n-2k$. We now can see that $e(f)(x)=y$, where $e(f)$ is the linear extension of $f$ to $\mathbb{Z}_2^n$ (see the proof of Theorem \ref{c2}).
 \\ (ii) Let $w(x) \neq w(y)$. Without loss of generality we can assume that $w(x)=2k-1$ and $w(y)=2k$. Let $x=e_{i_1}+...+e_{i_{2k-1}}$  and  $y=e_{j_1}+...+e_{j_{2k}}$. Note that
$y=(e_{j_1}+e_{j_{2k}})+(e_{j_2}+e_{j_{2k}})+...+(e_{j_{2k-2}}+e_{j_{2k}})+(e_{j_{2k-1}}+e_{j_{2k}}) $.
 There are vertices $e_{x_1},...,e_{x_{n-2k+1}}$ and
$e_{y_1}=e_{j_{2k}},e_{y_2},...,e_{y_{n-2k+1}}$ in $Q_n^2$ such that
$$\{ e_{i_1},...,e_{i_{2k-1}},e_{x_1},...,e_{x_{n-2k+1}}  \}=C_0, $$
$$\{e_{j_1}+e_{j_{2k}},e_{j_2}+e_{j_{2k}},...+e_{j_{2k-2}}+e_{j_{2k}},e_{j_{2k-1}}+e_{j_{2k}} \} \cup $$
$$\{ e_{y_1},e_{y_2}+e_{j_{2k}},...,e_{y_{n-2k+1}}+e_{j_{2k}} \}
 =C_{j_{2k}}$$
 We now   define the bijection $g$ from $C_0$ to $C_{j_{2k}}$ by the rule $g(e_{i_{r}}) =e_{j_{r}}+e_{j_{2k}}$, and $g(e_{x_1})=e_{y_1}$,  $g(e_{x_i})=e_{y_i}+e_{j_{2k}}$, $i \neq 1$. Let $e(g)$ be the linear extension of $g$ to $\mathbb{Z}_2^n$.  This yields that $e(g)$ is an automorphism of the graph $Q_n^2$ such that $e(g)(x)=y$.
\end{proof}

Theorem \ref{c3}  implies many results. For instance,
we now can deduce from it the following  corollary, which is   important in applied graph theory and interconnection networks.

%
%

\begin{cor} \label{c4} Let $n \geq 4$ be an integer. Then the connectivity of the graph $Q_n^2$ is maximal, namely,   $n$+$n\choose 2$ (its valency).

\end{cor}

\begin{proof}
 By Theorem \ref{c3} the graph $Q_n^2$ is distance-transitive, then it is edge-transitive. Thus, it follows from Theorem \ref{a},  that  the connectivity of the graph $Q_n^2$ is its valency, namely, $n$+$n\choose 2$.
\end{proof}
A block $B$, in the action of a group $G$ on a set $X$, is a subset of $X$ such that $B \cap g(B) \in \{B,\emptyset \}$,  for each $g$ in $G$. If $G$ is transitive on $X$, then we say that the permutation group $(X, G)$ is primitive if the only blocks are the trivial blocks, that is, those with cardinality 0,1 or $|X|$. In the case of an imprimitive permutation group $(X, G)$, the set $X$ is partitioned into a disjoint union of non-trivial blocks, which are permuted by $G$. We   refer to this partition as a block system.
A graph $\Gamma$ is said to be primitive or imprimitive according to   the group $Aut(\Gamma)$ acting on $V(\Gamma)$ has the corresponding property.  In the sequel, we need the following definition.
%
%

 \begin{defn}  \label{c5} A graph $\Gamma=(V,E)$ of diameter $D$ is said to be $antipodal$ if  for any    $u, v, w \in V$ such that $d(u, v) = d(u, w) = D$, then we have $d(v, w) = D$ or $v = w$. 

\end{defn}
Let $\Gamma_i(x)$ denote  the set of vertices of $\Gamma$ at distance $i$ from the  vertex $x$.                                                                                                                                        Let $\Gamma$ be a distance-transitive graph.  From Definition \ref{c5}  it follows that if $\Gamma_D(x)$ is a singleton set, then the graph $\Gamma$ is antipodal.
 It is easy to see that the hypercube  $Q_n$ is  antipodal, since every vertex $u$ has a unique vertex at maximum distance from it. Note that this graph is at the same time bipartite. We have the following fact \cite{1}.
 %
%
\begin{prop} \label{c6} A distance-transitive graph $\Gamma$ of diameter $D$ has a block $X = \{u\} \cup \Gamma_D(u)$ if and only if $\Gamma$ is antipodal, where $\Gamma_D(u)$ is the set of vertices of $\Gamma$ at distance $D$ from the vertex $u$.

\end{prop}

Also, we have the following important fact \cite{1}.
  %
%
\begin{thm} \label{c7} An imprimitive distance-transitive graph is either bipartite or antipodal. (Both possibilities can occur in the same graph.)

\end{thm}
We now can state and prove the following fact concerning the square of the hypercube $Q_n$.

 
\begin{cor} \label{c8} Let $n\geq 4$ be an integer. Then, the square of the hypercube $Q_n$, namely, the graph $Q_n^2$, is an imprimitive distance-transitive graph if and only if $n$ is an odd integer.

\end{cor}

\begin{proof} We know from Theorem \ref{c3}, that the graph $\Gamma=Q_n^2$ is a distance-transitive graph. Let $n=2k$ be an even integer. If $D$ denotes the diameter of $Q_n^2$, then $D=k$.  Let $C_0=\{e_1,...,e_n\}$ be the standard basis of the hypercube $Q_n$. Let $w=e_1+e_2+...+e_n$ and $B_1=\{w+e_i \ | \  1\leq i \leq n \}$. Consider the vertex $u=0$. It is easy to show that $\Gamma_D(u)=\{ w \}\cup B_1$. Two vertices $w$ and $w+e_1$ are in $\Gamma_D(u)$, but they are not at distance $k=D$ from  each other, since they are adjacent and $k>1$. Thus, when $n$ is an even integer, then the graph $Q_n^2$ is not antipodal. Since the girth of $Q_n^2$ is 3, then this graph is not bipartite. Now,   Theorem \ref{c7}  implies    that the graph $\Gamma=Q_n^2$ is not imprimitive. \

  Now assume that $n=2k+1$ is an odd integer. It is easy to see that $D=k+1$ and $\Gamma_D(0)=\{w\}$. Therefore by Proposition \ref{c6},  $\Gamma$ is antipodal, and hence has the set $\{0,w\}$ as a block. We now conclude that, when $n$ is an odd integer, then $Q_n^2$ is an imprimitive graph.
\end{proof}

Let $\Gamma=(V,E)$ be a simple connected graph with
diameter $D$.
A $distance$-$regular$  graph $\Gamma=(V,E)$, with diameter $D$, is a regular connected graph of valency $k$ with the following property. There are positive integers
$$b_0 = k, b_1, ...,b_{D-1};c_1=1, c_2,...,c_D, $$
such that for each pair $(u, v)$  of vertices satisfying $u \in \Gamma_i(v)$, we have\

(1) the number of vertices in $\Gamma_{i-1}(v)$ adjacent to $u$ is $c_i$, $1\leq i \leq D$. \

(2) the number of vertices in $\Gamma_{i+1}(v)$ adjacent to $u$ is $b_i$, $0\leq i \leq D-1$.\\
The intersection array of $\Gamma$ is $i(\Gamma) = \{k,b_1, ...,b_{D-1};  1,c_2, ...,c_d \}$.\

It is easy to show that if $\Gamma$ is a distance-transitive graph, then it is distance-regular \cite{1}. Hence, the hypercube $Q_n, \ n>2$ is a distance-regular graph.
We can verify  by an easy argument   that the intersection array of $Q_n$ is 
$$ \{n,n-1,n-2,...,1; 1,2,3,...,n  \}.$$
In other words, for hypercube $Q_n$, we have $b_i=n-i$,   $c_i=i$, $1 \leq i \leq n-1$, and $b_0=n$, $c_n=n$. In the following theorem, we determine the intersection array of the square of the hypercube $Q_n$ \cite{1}.
%
%

\begin{thm} \label{c9} Let $n>3$ be an integer and $\Gamma=Q_n^2$ be the square of the hypercube $Q_n$. Let    $D$ denote  the diameter of $Q_n^2$. Then for the intersection array of this graph  we have $b_0$=$n+1 \choose 2$,  $b_i$=$ n-2i+1 \choose 2$, $c_i$=$2i \choose 2$, $1\leq i \leq D-1$.
 Also,  $c_D$=$n+1 \choose 2$, when $n$ is an odd integer and $c_D$=$n \choose 2$ when $n$ is an even integer.

\end{thm}

\begin{proof}
Since   $Q_n^2$ is a regular graph of valency $n+1 \choose 2$, thus we have
$b_0$=$n+1 \choose 2$. Let $u$ be a vertex in $Q_n^2$ at distance $i$ from the vertex $v=0$. It is easy to check that $w(u)=2i$ or $w(u)=2i-1$. This implies  that  that the diameter of the graph
$Q_n^2$ is $D=  \lceil \frac{n}{2} \rceil$.\

Let $u$ be a vertex in $Q_n^2$ at distance $i \geq 1$ from the vertex $v=0$, such that $i \neq  D$. There are two cases, that is,  $w(u)=2i$,  or  $w(u)=2i-1$. Without lose of generality we can assume that $w(u)=2i$.
 Hence $u$ is of the form $u=e_{j_1}+e_{j_2}+...+e_{j_{_{2i}}}$. Now it is easy to show that if $x$ is a vertex of $Q_n^2$ adjacent to $u$  and at distance $i-1$
from  the vertex $v=0$, then $x$ must be of the
form $x=u+e_k+e_l$, where $e_k,e_l \in  \{ e_{j_1},e_{j_2},...,e_{j_{_{2i}}} \}$. It is clear that the number of such $xs$ is equal to    $2i \choose 2$. Moreover, If $x$ is a vertex of $Q_n^2$ adjacent to $u$  and at distance $i+1$
from  the vertex $v=0$, then $x$ must be of the
 forms $x=u+e_k$ or $x=e_k+e_l$, where $e_k,e_l \in \{ e_1,e_2,...,e_n  \} - \{ e_{j_1},e_{j_2},...,e_{j_{_{2i}}} \}$. It is clear that the number of such $xs$ is equal to    $n-2i \choose 1$+$n-2i \choose 2$=$n-2i +1 \choose 2$.
We now deduce that when $1\leq i \leq D-1$, then
$c_i$=$2i \choose 2$, and $b_i$=$n-2i+1 \choose 2$. \

When $n$ is an odd integer, then the vertex $u=e_1+e_2+...+e_n$ is the unique vertex of $Q_n^2$ at distance
$D$ from the vertex $v=0$. Thus $c_D$=$n+1 \choose 2$, namely,  the valency of $u$. If $n$ is
an even integer, then $\Gamma_D(0)=\{u,u+e_i | \ 1 \leq i \leq n  \}$ is the set of vertices of $\Gamma=Q_n^2$ at distance $D$ from the vertex $v=0$. Now, by a similar argument which is done in the first section of the proof, it can be shown that  $c_D$=$n \choose 2$.
\end{proof}

\begin{rem} \label{c10} There are   distance-regular graphs $\Gamma=(V,E)$, with the property  that their squares are not  distance-regular. For instance, consider the cycle $C_n$ with vertex set $\{0,1,2,...,n-1  \}$. It is well known that $C_n$ is a distance-regular graph of diameter $[ \frac{n}{2} ]$ with the intersection array:  \

$\{2,1,1,...,1,1;1,1,1,...,1,2\}$ when $n$ is an even integer   and, \
 
$\{2,1,1,...,1,1;1,1,1,...,1,1\}$ when $n$ is an odd integer \cite{1}. \\ Now,  assume that $n\geq 7$. It can be shown by an easy argument that
$\Gamma=C_n^2$ is not a distance-regular graph. To see this fact, let $v$ be a vertex in $C_n$ at distance $i$ from the vertex $0$, and  $c_i(v)= |\Gamma_{i-1}(0) \cap N(v)|$.  It is easy to show that  $\Gamma_i(0)=\{2i, -2i, 2i-1, -2i+1  \}$, and  $c_i(2i)=1$, but $c_i(2i-1)=2$.
\end{rem}

\begin{rem} \label{c11}
Let  $n,k \in \mathbb{ N}$ with $ k < n,   $ and let $[n]=\{1,...,n\}$.  Consider the  Johnson  graph  $J(n,k)$. It is clear that the order of this graph is $n \choose k$.  It is easy to check that $J(n,k) \cong J(n,n-k)$, hence we assume that $1 \leq k \leq \frac{n}{2}$.
The class of  Johnson graphs is one of the most well known and interesting subclass of   distance-regular graphs \cite{3}. It is easy to show that if $v$ and $w$ are vertices in the Johnson graph $J(n,k)$, then $d(v,w)=k-|v \cap w|$. Thus, the diameter of the Johnson graph $J(n,k)$ is $k$. Note that the graph $J(n,1)$ is the complete graph $K_n$ and hence it is distance-regular. The diameter of the graph $J(n,2)$ is 2, hence the diameter of its square is 1. Thus the graph $J^2(n,2)$ is the complete graph $K_m$, and hence it is a distance-regular graph ($m$=$n \choose 2$).
We can show that when $k=3$, then the square of Johnson graph $\Gamma=J(n,k)$ is a distance-regular graph if and only if $n=6$. For checking this,
let $v=\{1,2,3\}$. Note that the  diameter of the
graph $\Gamma^2$ is 2 and a vertex $w$ in
$\Gamma^2$ is at distance 2 from $v$ if and only if $|v \cap w|=0$. Moreover $w$ is at distance 1 from   $v$ if and only if
$|v\cap w|\in \{ 1,2 \}$. Hence $\Gamma^2_1(v)=V(\Gamma)-\{v,v^c\}$  and $\Gamma^2_2(v)=\{v^c\}$, where $v^c$ is the complement of the set $v$ in the set $\{1,2,...,6\}$. Thus  $v^c=\{4,5,6\}$. Now, it is clear that $b_0(v)$=${3}\choose {2}$${3}\choose {1}$+${3}\choose {1}$${3}\choose {2}$=18. Also,  for every $w \in \Gamma^2_1(v)$,  $c_1(w)=1$ and $b_1(w)=1$, and   $c_2(v^c)=|\Gamma^2_1(v)|=18$. Thus the graph $\Gamma^2 =J^2(6,3)$ is a distance-regular
graph with intersection array $\{18,1;1,18\}$. But, if $n>6$, then the graph $\Gamma^2= J^2(n,3)$ is not distance-regular. In fact if
$n>6$, then for the vertex $v=\{1,2,3\}$,
each of the vertices $u=\{1,2,4\}$ and $w=\{ 1,4,5\}$ is in $\Gamma^2_1(v)$. If $x \in \Gamma^2_2(v)$ is adjacent to $u$, then $4 \in x$, and hence $x=\{4\}\cup y$, where $y \subset v^c- \{4\}$ and $|y|=2$. We now can deduce that  $b_1(u)$=${n-4}\choose {2}$. On the other hand,  if $x \in \Gamma^2_2(v)$ is adjacent to $w$, then $4 \in x$ and $5 \notin x$, or $5 \in x$ and $4 \notin x$ or $4,5 \in x$. Thus, $b_1(w)$=2${n-5}\choose {2}$+${n-5} \choose {1}$=${n-4}\choose {2}$+${n-5}\choose {2}$. This implies that when $n \geq 7$  then the graph $J^2(n,3)$ cannot be distance-regular.\

By a similar argument  we  we can show that the graph $J^2(8,4)$ is distance-regular, but if $n >8$, then the graph $J^2(n,4)$ is not distance-regular.
\end{rem}

\begin{rem} \label{c12}  Let $\Gamma=(V,E)$ be a graph.
$\Gamma$ is said to be a $strongly$ $regular$ graph with parameters
 $(n,k,\lambda,\mu)$, whenever $|V|=n$, $\Gamma$
is a regular graph of valency $k$, every pair of adjacent vertices of $\Gamma$ have $\lambda$ common neighbor(s), and every pair of non-adjacent vertices of $\Gamma$ have $\mu$ common neighbor(s). It is clear that the diameter of every strongly regular graph is 2. It is easy to show that if a graph $\Gamma$ is a distance-regular graph of diameter 2 and order $n$, with intersection array $(b_0,b_1;c_1,c_2)$, then $\Gamma$ is a strongly regular graph with parameters
$(n,b_0,b_0-b_1-1,c_2).$ We know that the diameter of the graph $Q_n^2$ is $\lceil \frac{n}{2} \rceil$. Now, it follows from Theorem \ref{c3}, that $Q_3^2$ is a strongly regular graph with parameter $(8,6,4,6)$. This graph is known as the $coktail$-$party$ graph $CP$(4) \cite{1}. Also, the graph $Q_4^2$ is a strongly regular graph with parameter $(16,10,6,6)$. We know that when a graph $\Gamma$ is is a strongly regular graph with parameters
 $(n,k,\lambda,\mu)$, then its complement is again a strongly regular graph with parameter $(n,n-k-1,n-2-2k+\mu,n-2k+\lambda)$ \cite{9}. Hence,  the complement of the graph $Q_4^2$ is a strongly regular graph with parameter $(16,5,0,2)$. This graph
 is known as the $Clebsch$ graph \cite{9} and it is the unique strongly regular graph with parameters $(16,5,0,2)$. Figure 1. displays a version of  the Clebsch graph (the complement of the graph $Q_4^2$) in the plane \cite{9}.
\end{rem}

\begin{figure} [ht]
\centerline{\includegraphics[width=11 cm]{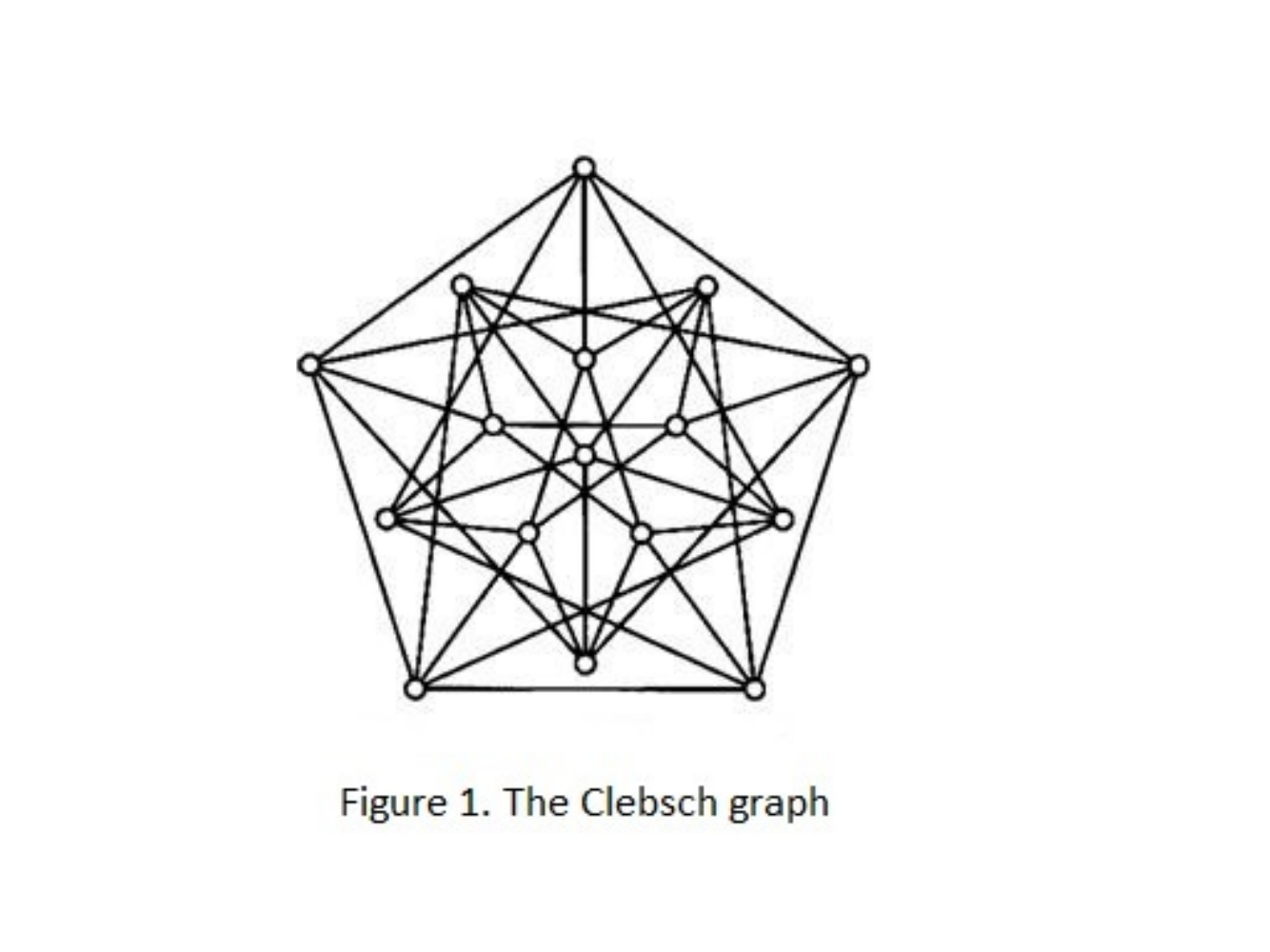}} 
\end{figure}

\section{The spectrum of the square of the hypercube}
The square of the hypercube $Q_n$ has some further  interesting algebraic properties. For obtaining some of those properties,  we need the spectrum of  this graph.  The spectrum of $Q_n$ is known \cite{1}, however we are not aware of a paper showing the spectrum of $Q_n^2$. Here we compute by means of an algebraic and self-contained method  the spectrum of $Q_n^2$.    Hence,  in the sequel and in the first step,  we determine  the spectrum of the graph $Q_n^2$.\

Let $\Gamma=(V,E)$ be a graph with the vertex set $\{v_{1},\cdots,v_{n}\}$. Then the adjacency matrix of $\Gamma$ is an $n\times n$ matrix $A=(a_{ij})$, in which  columns and rows are labeled by $V$ and $a_{ij}$ is defined as follow:
\begin{center}
$a_{ij}=A(v_i,v_j)=\begin{cases}1 &  $if$ \  v_i \    $is$  \  $adjacent$ \    $to$  \    v_j    \\ 0 &  $otherwise.$  \end{cases} $
\end{center}
 If $Ax = \lambda x, x \neq 0$,
then $\lambda$  is an eigenvalue of $A$, and $ x$ is an eigenvector of $A$ corresponding to $\lambda$ \cite{9}.   Let $\lambda_{1},\cdots,\lambda_{r}$ be eigenvalues of $A$ with multiplicities $m_{1},\cdots,m_{r}$, respectively. The spectrum of the  graph $\Gamma$ is defined   as 
\begin{center}
$Spec(\Gamma)=\left\lbrace\begin{array}{cccc}\lambda_{1}&\lambda_{2}&\cdots&\lambda_{r}\\
m_{1}&m_{2}&\cdots&m_{r}
\end{array}\right\rbrace. $
\end{center}

When we work with graphs there is an additional refinement. We can suppose that an eigenvector is a  real function $f$  on the vertices. Then if at any vertex $v$  you
sum up the values of $f$ on its neighboring vertices, you should get $\lambda$  times the values of $f$ at
$v$.  Formally,  \

$$\sum_{w \in N(v)} {f(w)} = \lambda f(v). $$ \

Let $G$  be a finite abelian group (written additively) of order $|G|$ with
identity element  0=$0_G$.  A  character $\chi$  of $G$  is a homomorphism from $G$ into the
multiplicative group $U$ of complex numbers of absolute value 1, that is, a
mapping from $G$ into $U$ with $\chi(g_1+g_2)= \chi(g_1) \chi(g_2)$ for all $g_1,g_2 \in G$.  If $G$ is a finite abelian group, then there are integers $n_{1},\cdots ,n_{k}$, such that $G=\mathbb{Z}_{n_{1}}\times \cdots \times \mathbb{Z}_{n_{k}}$.
 Let   $S=\{s_{1},\cdots ,s_{n}\} $ be a non-empty subset of $G $ such that $0\not\in S $ and $S=-S$. Let $\Gamma =Cay(G;S)$. Assume $f:G\longrightarrow \mathbb{C}^{*}$ is a character where $\mathbb{C}^{*}$ is the multiplicative group of the complex numbers.
 If $\omega _{ij}=e^{\frac{2\pi ij}{n_{i}}},\; 0\leq i\leq k,  \; 1\leq j\leq n_{i}$, is an $n_{i}$th root of unitary, then  $f$   is of the form
$f=f_{(\omega _{1},\cdots ,\omega _{k})}$, where  $f_{(\omega _{1},\cdots ,\omega _{k})}(x_{1},\cdots ,x_{k})=\omega _{1}^{x_{1}}\omega _{2}^{x_{2}}\cdots \omega _{k}^{x_{k}}$, for each $(x_1,x_2,...,x_k) \in G$ \cite{11}.

       If $v $ is a vertex of $\Gamma $, then we know that $N(v)=\{v+s_{1},\cdots , v+s_{n}\}$ is the set of vertices that are adjacent to $v$. We now have

$$\sum _{w\in N(v)} f(w)=\sum _{i=1}^{n}f(v+s_{i})=\sum_{i=1}^{n}f(v)f(s_{i})=f(v) (\sum_{i=1}^{n}f(s_{i}) ).$$

Therefore, if we let $\lambda =\lambda _{f}=\sum _{s\in S}f(s)$ then we have $\sum _{w\in N(v)}f(w)=\lambda _{f}f(v)$, and hence the mapping $f$ is an eigenvector for the Cayley graph $\Gamma $ with corresponding eigenvalue $\lambda =\lambda _{f} =\sum _{s\in S}f(s)$.
%
%

\begin{thm} \label{d1}
Let $n >3$ be an integer and $Q_n^2$ be the square of the hypercube $Q_n$. Then each of the eigenvalues  of  $Q_n^2$ is of the form,
 $$\lambda_i=\frac{1}{2}n(n+1)-2i(n+1)+2i^2, $$
  for $0 \leq i \leq
\lfloor \frac{n+1}{2} \rfloor$. Moreover,    the multiplicity of $\lambda_0$ is 1, the multiplicity of $\lambda_i$ is
$m(\lambda_i)$= $n \choose i $+$n \choose n+1-i$,
 for $1 \leq i \leq
\lfloor \frac{n+1}{2} \rfloor$, when $n$  is an  even integer,
 and  $m(\lambda_i)$=$n \choose i $+$n \choose n+1-i$  for $1 \leq i <
\lfloor \frac{n+1}{2} \rfloor$, when $n$  is  an odd integer,  with $m(\lambda_j)$=$n \choose j$ for $ j=\lfloor \frac{n+1}{2} \rfloor$.

\end{thm}
\begin{proof}
According to what is stated before this theorem, every eigenvector of the graph $\Gamma=Q_n^2=Cay(\mathbb{Z}_{2}^{n};S)$ is of the form $f=f_{(\omega_{1}, \cdots ,\omega_{n})}$, where each $\omega_{i} $, $1\leq i\leq n$, is a complex number such that $\omega _{i}^{2}=1$, namely, $\omega_{i}\in \{1,-1\}$.  We now have 
$$\lambda_f= \sum_{w \in S} {f(w)} = \sum_{i=1}^n {f(e_i)}  +\sum_{i,j=1, \  i \neq j}^n  f(e_i+e_j)$$
$$= \sum_{i=1}^n {f(e_i)}   +\sum_{i,j=1, \   i\neq j}^n {f(e_i)f(e_j)}.$$
Note that
  for every vertex $v=(x_1,\ldots,x_n)$,  $x_i\in\{0,1\}$ in $Q_n^2$, we have
$$f(x_1,\ldots,x_n)
=f_{(w_1,\ldots,w_n)}(x_1\ldots,x_n)=w_1^{x_1}\ldots w_n^{x_n}.$$
 Note that  in the computing of the value of $w_1^{x_1}\ldots w_n^{x_n}$ we can ignore $w_i$ when $w_i=1$.   Thus,   for $e_k=(0,\ldots,0,1,0\ldots,0)$, where $1$ is the $k$th entry,   we   have; \

$$f(e_k)=f_{(w_1,\ldots,w_n)}(0,\ldots,0,1,0,\ldots,0)$$
$$=w_1^0\ldots w_k^1w_{k+1}^0\ldots w_n^0=
\begin{cases}
 \ -1 \  \  if\  \  w_k=-1\\
 \  \ 1\  \ \

\
 if \   \ w_k=1
\end{cases}$$
 \\
Hence, if in the $n$-tuple $(w_1,\ldots,w_n)$ the number of $-1$$s$ is $i$ (and therefore the number of  l$s$ is $(n-i)$), then in the sum
$$\sum_{k=1}^{n}f(e_k)=\sum_{k=1}^{n}f_{(w_1,\ldots,w_n)}(0,\ldots,x_k,0,\ldots,0), \ \ x_k=1, $$
 the contribution of  $-1$ is $i$  and the contribution of 1 is $n-i$. Therefore,   we have
  $$\sum_{k=1}^{n}f(e_k)=-i+(n-i)=n-2i.$$
On the other hand,  since 

$$(\sum_{k=1}^{n}f(e_k))^2=\sum_{k=1}^{n}{f(e_k)}^2+2\sum_{i,j=1,\    i\neq j}^n {f(e_i)f(e_j)}, $$
 therefore,  we have 

$$\sum_{i,j=1,\    i\neq j}^n {f(e_i)f(e_j)}=\frac{1}{2}({(n-2i)}^2- \sum_{k=1}^{n}{f(e_k)}^2). $$
 Now since $\sum_{k=1}^{n}{f(e_k)}^2=n$, thus we have 
  $$ \lambda_f=\sum_{i=1}^n {f(e_i)}+\sum_{i,j=1, \   i\neq j}^n {f(e_i)f(e_j)}
={(n-2i)}+\frac{1}{2}({(n-2i)}^2-n)$$
$$=\frac{1}{2}n+\frac{1}{2}n^2-2ni
+2i^2-2i=\frac{1}{2}n(n+1)-2i(n+1)+2i^2.$$
Note that $f=f_{(w_1,w_2,...,w_n)}$, and the number of sequences $(w_1\ldots,w_n)$   in which $i$ entries are $-1$ is ${n\choose i}$.
If we denote $\lambda_f$ by $\lambda_i$, then we deduce that every eigenvalue of the graph $Q_n^2$ is of the form 
$$\lambda_i= \frac{1}{2}n(n+1)-2i(n+1)+2i^2, \ 0 \leq i \leq n. \ \ \ (**)$$
Consider the real function $f(x)=\frac{1}{2}n(n+1)-2x(n+1)+2x^2$. Then $\lambda_i=f(i)$, $i \in \{0,1,...,n \}$. This function reaches its minimum at $x=\frac{n+1}{2}$. Now by using some calculus,  we can see that $f(x)=f(n+1-x)$. Thus, we have $\lambda_i=f(i)=f(n+1-i)=\lambda_{n+1-i}$, $1 \leq i \leq n$. Now  it  follows that if $n=2k$, then the multiplicity of $\lambda_i$ is $n \choose i$+$n \choose n+1-i$,   $1 \leq i \leq k$. Note that when $n=2k+1$, then $n+1-(k+1)=k+1$, thus $\lambda_{n+1-(k+1)}=\lambda_{k+1}$.
Hence  if $n=2k+1$,  then the multiplicity of $\lambda_i$ is $n \choose i$+$n \choose n+1-i$,   $1 \leq i \leq k$, and the multiplicity of $\lambda_{k+1}$ is $n \choose k+1$. Note that since the graph $Q_n^2$ is a $n+1 \choose 2$-regular graph,  hence the multiplicity of $\lambda_0$=$n+1 \choose 2$=$\frac{1}{2}(n+1)n$ is 1.
\end{proof}

Let $ \Gamma =(V,E)$ be a graph.  The line graph $ L( \Gamma )$
of the graph $ \Gamma$ is constructed by taking the edges of $ \Gamma$
as vertices of  $ L( \Gamma )$, and joining two vertices in $ L( \Gamma )$
whenever the corresponding edges in $ \Gamma $ have a common vertex. Note that if $e=\{ v,w \}$ is an edge of $\Gamma$, then its degree in the graph $L(\Gamma)$ is $deg(v)+deg(w)-2$.
Concerning the eigenvalues of the  line graphs, we have the following fact \cite{1,9}.

\begin{prop} \label{d2}
 If $\lambda$ is an eigenvalue of a line graph $L(\Gamma)$, then $\lambda \geq -2$.

\end{prop}

Therefore, if $\lambda < -2$  is an  eigenvalue of a  graph graph $\Gamma$,  then $\Gamma$ is not a line graph. \

A $(c, d)$-$biregular$  graph is a bipartite graph in
which each vertex in one part has degree $c$ and each vertex in the other part has degree $d$ \cite{25}. It is known and easy to prove that if the line graph of the graph $\Gamma$ is regular, then $\Gamma$ is a regular or a  $(c, d)$-biregular  bipartite graph.

\begin{thm} \label{d3} Let $n\geq4$ be an integer and
 $Q_n^2$ be the square of the hypercube $Q_n$. Then $Q_n^2$ cannot be  a line graph.

\end{thm}

\begin{proof}
Let $k=\lfloor \frac{n}{2} \rfloor$. Hence, if $n$ is an
even integer, then $n=2k$ and if $n$ is an odd integer then $n=2k+1$. It follows from Theorem \ref{d1}, that the smallest eigenvalue of the graph $Q_n^2$ is  $
\lambda_k$, when $n$ is an even integer and $\lambda_{k+1}$, when $n$ is an odd integer. Now
  consider the eigenvalue $\lambda_k$ of the graph $Q_n^2$  in (**) (in the proof of Theorem \ref{d1}). Therefore if $n$ is an even integer, then
we have 

$$\lambda_k=k(2k+1)-2k(2k+1)+2k^2=k(2k+1-4k-2+2k)=-k.$$

Moreover if $n=2k+1$, then  we have,

$$\lambda_{k+1}=(2k+1)(k+1)-2(k+1)(2k+2)+2{(k+1)}^2=$$

$$(k+1)(2k+1-4k-4+2k+2)=-k-1.$$
\

We now deduce that when  $n\geq 5$, then $\lambda_k \leq -3$.
Now, it follows from Proposition \ref{d2},  when $n\geq 5$, then the graph $Q_n^2$ can not be a line graph. \

Our argument shows that if $\lambda$ is an eigenvalue of the graph $Q_4^2$, then $\lambda \geq -2$, and hence in this way we can not say anything about our claim. \\We now show that $Q_4^2$ is not a line graph. On the contrary, assume that $Q_4^2$ is a line graph. Thus, there is a graph $\Delta$ such that $Q_4^2=L(\Delta)$. Since $Q_4^2$ is a regular graph, hence  (i) $\Delta$ is a  regular graph,  or (ii)  $\Delta$ is  a  biregular bipartite graph.\\
(i) Let $\Delta=(V,E)$ be a $t$-regular graph of order $h$. Since $Q_4^2$ is 10-regular, thus, $L(\Delta)=Q_4^2$ is a $2t-2=10$-regular graph, and hence $t=6$. Therefore we have $16=|E|
=\frac{1}{2}6h=3h$, which is impossible.\\ (ii) Let
$\Delta=(A\cup B,E)$ be a $(c,d)$-biregular bipartite   graph such that every vertex in $A$ $( B )$ is of degree $c$ $( d )$.   Hence we have $16=|E|=c|A|=d|B|$. Thus $c$ and $d$ must divide 16. On the other hand, if $e=\{a,b \}$ is an edge of $\Delta$, then we must  have $deg(a)+deg(b)-2=10=c+d-2$. Hence we have $c+d=12$. We now can check that $\{c,d\}=\{4,8\}.$ Without loss  of generality, we can assume that $d=8$ and $c=4$. Hence we must have $|A|\geq 8$. Now since each vertex in $A$ is of degree $c=4$, then we must   have, $16 =|E|=c|A|=4|A| \geq 4 \times 8$=32, which is impossible.  \

Our argument shows that the graph $Q_4^2$ is also not  a line graph.
\end{proof}

An $automorphic$ graph is a distance-transitive graph whose automorphism group acts primitively on its vertices, and not a complete graph or a line graph.\\
Automorphic graphs are apparently very rare. For instance, there are exactly   three cubic automorphic graphs \cite{1,2}.  It is clear that  for  $n\geq 3$,    the graph $Q_n^2$ is not a complete graph. We now derive from Corollary \ref{c8}, and Theorem \ref{d3},   the following important result.

\begin{cor} \label{d4} Let $n\geq 4$ be an integer. Then the square of the hypercube $Q_n$, that is,  the graph $Q_n^2$,  is an automorphic graph if and only if $n$ is an even integer.
\end{cor}

\section{Conclusion} In this paper, we proved that the square of the distance-transitive graph $Q_n$, that is,  the graph $Q_n^2$,  is again a distance-transitive graph (Theorem \ref{c3}). We  showed that there are important classes of  distance-transitive graphs (including the cycle $C_n$, $n \geq 7$),   such that their squares are not even
distance-regular  (and hence are not distance-transitive) (Remark \ref{c11}). Also, we determined the spectrum of the graph $Q_n^2$ (Theorem \ref{d1}). Moreover, we   showed that  when $n>3$ is an even integer, then the graph $Q_n^2$ is an automorphic graph, that is, a distance-transitive primitive graph which is not a complete or a  line graph (Corollary \ref{d4}).


\begin{thebibliography}{}
\bibitem{1}  N. L. Biggs,  Algebraic graph theory. 2nd ed. Cambridge: Cambridge Mathematical Library; 1993.
(Cambridge University Press).
\bibitem{2}N. L. Biggs,  D.H. Smith,   On trivalent graphs, Bull. London Math, Soc. 3,  (1971),  155-158.



\bibitem{3} A. E.  Brouwer,  A. M.  Cohen,  A. Neumaier, Distance-Regular Graphs, Springer-
Verlag, New York, (1989).
\bibitem{4}G. Chartrand, A. M. Hobbs,   H. A. Jung,  S. F. Kapoor,    J. A. Nash-Williams, The square of a block is Hamiltonian-connected. J. Combin. Theory
Ser. B 16 (1974),  290-292.

\bibitem{5} R. Diestel,  Graph Theory (4th ed.), Springer-Verlage,  Heildelberg (2010).


 \bibitem{6}  J. D.  Dixon,   B. Mortimer,   Permutation Groups,  Graduate Texts in Mathematics 1996; 163:   Springer-Verlag, New
York.



\bibitem{7} H. Fleischner, In the square of graphs, Hamiltonicity  and pancyclicity, hamiltonian connectedness and panconnectedness are equivalent concept, Monatsh.
Math. 62 (1976),  125-149.

\bibitem{8}A. Ganesan,  Automorphism group of the complete transposition graph. J. Algebr. Comb. (2015).
https://doi.org/10.1007/s10801-015-0602-5.

\bibitem{9}  C. Godsil,    G. Royle,   Algebraic Graph Theory,   Berlin: Springer  2001.


\bibitem{10} X. Huang,   Q. Huang,   Automorphism group of the complete alternating group, Appl. Math. Comput. 314 (2017), 58-64.

\bibitem{11} G. James, M. Liebeck,  Representations and Characters of Groups, Cambridge University Press (2001).

\bibitem{12} G. A. Jones, R. Jajcay, Cayley properties of merged Johnson graphs. J. Algebr. Comb. 44, (2016),  1047–1067.







\bibitem{13} S. M. Mirafzal,    Some other algebraic properties of folded hypercubes,  Ars Comb. 124   (2016), 153-159.




\bibitem{14} S. M. Mirafzal,  More odd graph theory from another point of view. Discrete Math. 341 (2018), 217–220.




\bibitem{15}S. M.  Mirafzal,   A new class of  integral  graphs constructed from the hypercube.    Linear Algebra Appl. 558 (2018),  186-194.





\bibitem{16} S. M.  Mirafzal,   The automorphism group of the bipartite Kneser graph, Proc. Math. Sci (2019), doi.org/10.1007/s12044-019-0477-9.

\bibitem{17}S. M. Mirafzal,  A. Heidari,  Johnson graphs are panconnected, Proc. Math. Sci (2019),  doi.org/10.1007/s12044-019-0527-3.

\bibitem{18}S. M.  Mirafzal,  M. Ziaee,  Some algebraic aspects  of enhanced Johnson graphs, Acta Math. Univ. Comenianae. 88(2) (2019), 257-266.

\bibitem{19}S. M. Mirafzal, Cayley properties of the line graphs induced by consecutive layers of the hypercube,
Bull. Malaysian Math. Sci. (2020), https://doi.org/10.1007/s40840-020-01009-3.


\bibitem{20}S. M. Mirafzal,  On the automorphism groups of connected bipartite irreducible graphs. Proc. Math.
Sci  (2020). https://doi.org/10.1007/s12044-020-0589-1.


\bibitem{21}S. M.  Mirafzal,  M. Ziaee, A note on the automorphism group of the Hamming graph,
Trans. Comb. 10 (2),  (2021), 129-136.

\bibitem{22}S. M. Mirafzal,  On the automorphism groups of us-Cayley graphs, arXiv: 1910.12563.v4 1702.02568v4 [math.GR] 11 May 2021.

\bibitem{23} S. M.  Mirafzal,   A note on  the automorphism groups of  Johnson graphs, Ars Comb. 154  (2021),   245-255 (Available from: arXiv: 1702.02568v4).


\bibitem{24}S. M. Mirafzal,   The automorphism group of the Andr\'{a}sfai  graph, Discrete Math. Lett. 10 (2022) 60-63.


\bibitem{25} E. R.  Scheinerman,  D. H.  Ullman,   Fractional graph theory, Wiley-Interscience Series in Discrete Mathematics and Optimization, New York: John Wiley and  Sons Inc, (1997).

\bibitem{26}Y. Wang, Y. Q.  Feng, Half-arc-transitive graphs of prime-cube order
of small valencies, Ars Math. Contemp 13  (2017),   343-353.

\bibitem{27} Watkins  M,  Connectivity of transitive graphs,  J. Combin. Theory 8 (1970), 
23-29.

\bibitem{28}J. X. Zhou,  J. H. Kwak,   Y. Q. Feng,  Z.L. Wu,   Automorphism group of the balanced hypercube, Ars Math. Contemp 12  (2017),   145-154.  

\end{thebibliography}
\end{document}